\documentclass[a4paper, twoside, notitlepage, 12pt]{amsart} 
\usepackage{amsmath,amsfonts,amssymb,amsthm,graphicx,mathrsfs} 

\usepackage{amscd}
\usepackage{xypic}
\xyoption{all}

\usepackage{verbatim}
\usepackage[latin1]{inputenc}   

\title[Non-effective Deformations of the Hilbert Functor]{Non-effective Deformations of Grothendieck's Hilbert Functor}
\author{Christian Lundkvist}
\address{Department of Mathematics \\ KTH \\ Stockholm \\ Sweden}
\email{chrislun@math.kth.se}

\author{Roy Skjelnes}
\address{Department of Mathematics \\ KTH \\ Stockholm \\ Sweden}
\email{skjelnes@math.kth.se}

\subjclass[2000]{Primary 14C05; Secondary 14B12, 14A20}

\DeclareMathOperator{\Spec}{Spec}

\DeclareMathOperator{\can}{can}

\newcommand{\ra}{\longrightarrow}
\newcommand{\sra}{\rightarrow}

\newcommand{\mas}[1]{ \left\{ #1 \right\} }

\newcommand{\del}{\delta}

\newcommand{\eqbeg}{\begin{equation}}
\newcommand{\eqend}{\end{equation}}
\newcommand{\arbeg}[1]{\begin{array}{#1}}
\newcommand{\arend}{\end{array}}

\newcommand{\m}{\mathfrak{m}}

\newcommand{\bimap}{\ensuremath{\rightrightarrows}}

\newtheorem{thm}[subsection]{Theorem}
\newtheorem{lemma}[subsection]{Lemma}
\newtheorem{cor}[subsection]{Corollary}
\newtheorem{prop}[subsection]{Proposition}

\theoremstyle{definition}

\newtheorem{exam}[subsection]{Example}

\theoremstyle{remark}
\newtheorem{rem}[subsection]{Remark}

\numberwithin{equation}{subsection}

\newcommand{\calH}{\mathcal{H}ilb}
\newcommand{\calS}{\mathcal{S}}
\newcommand{\calO}{\mathcal{O}}
\newcommand{\calF}{\mathcal{F}}

\begin{document}
\pagenumbering{arabic}

\begin{abstract} 
We show that the Hilbert functor of rank one families on a non-separated scheme $X$ admits deformations that are not effective. For such ambient schemes  we have that the Hilbert functor is not representable by a scheme or an algebraic space.
\end{abstract}

\maketitle


\section*{Introduction}

One of Artin's criteria for representability of functors is the condition of effectivity of formal deformations \cite{artin_alg_form_moduli_I}. We will in this note show that the Hilbert functor of Grothendieck \cite{fga} does not always fulfill this criterion, and in particular that those functors will not be representable by an algebraic space, or a scheme.

It has been implicitly known for some time that the Hilbert functor parametrizing flat families does not behave well when the ambient scheme $X$ is not separated. In this note we show explicitly that separatedness is a neccessary condition for the Hilbert functor to be representable.

In fact, when $X$ is not separated over the base $S$ we show that the subfunctor $\calH^1_{X/S}$ parametrizing rank one families is not representable. In the separated case the functor is represented by the scheme $X$ itself, as intuition would suggest.\\
\\
We briefly sketch the arguments in the paper. Fix the ambient scheme $X$, in the category of schemes over a fixed base scheme,  and we let $\calH^1_{X}$ denote the Hilbert functor of 1-points on $X$, suppressing the base scheme in the notation. That is, the functor $\calH^1_{X}$ parametrizes families of closed subschemes of $X$ that are flat and of finite rank one over the base. The crucial fact in the definition is that the families are closed. 

Furthermore, for a complete local ring $A$ we consider the natural map
\[
\calH^1_X(A) \ra \varprojlim \calH^1_X(A/\m^{n+1}). 
\]
A formal deformation $\{\xi_n\}$ is a collection of compatible families $\xi_n\in \calH^1_X(A/\m^{n+1})$, and the deformation is called effective if it is in the image of the above map. It is easy to see that the map above is injective, but the map is not surjective in general. When the fixed scheme $X$ does not satisfy the valuative criterion of separatedness we show that there exists a complete valuation ring $A$ such that surjectivity of the above map fails.  In particular we have that the Hilbert functor $\calH_X^1$ is not representable for such schemes.

An explanation for the above mentioned result is as follows. It is easy to see that $\calH^1_X$ parametrizes \emph{closed} sections of the structure map $f : X \sra S$ and when the morphism $f$ is separated we have that any section is closed. On the other hand, schemes $f : X\sra S$ that do not satisfy the valuative criterion of separatedness have non-closed sections. Replacing $S$ with the spectrum of a complete valuation ring $A$, we have that different extensions of the generic point of the curve $S$ yield sections $\xi : S\sra X$ that are not closed. However, the infinitesimal truncations $\xi_n : \Spec (A/\m^{n+1}) \sra X$ of a section $\xi$ are closed. Consequently the infinitesimal truncations $\xi_n$ form a formal deformation, which is not effective since $\xi \notin \calH^1_X(A)$ and the section $\xi$ is uniquely determined by the $\{\xi_n\}$.

In the following sections we give the details of the above sketch.


\section{The Hilbert functor of one point}

We will in this first section define Grothendieck's Hilbert functor of points, and recall what is known for the Hilbert functor of one point.

\subsection{Sections} 
\label{sec:sections}
We fix  a morphism of schemes $f : X \sra S $. A \emph{section} $\xi $ is a morphism of schemes $\xi : S \sra X$ such that the composition $f\circ \xi$ is the identity on $S$. Sections are always immersions \cite[Cor. 5.3.11]{egaI}, and we say that a section is \emph{closed} if it is a closed immersion.

We let $\calS_X$ denote the contravariant functor that to any $S$-scheme $T$ assigns the set $\calS_X(T)$ of sections $T\sra X \times_S T$ of the projection map $X\times_S T\sra T$.

\begin{lemma}
\label{lem:sec_to_hom}
Composing a section $T \sra X \times_S T$ with the projection map $X\times_ST \sra X$ gives a map of functors
\eqbeg
\label{eq:sec_to_hom}
n: \calS_X \ra \mathrm{Hom}_S(-,X),
\eqend
which is an isomorphism. The inverse of (\ref{eq:sec_to_hom}) is induced by the diagonal map $\Delta : X \sra X\times_S X$.
\end{lemma}

\begin{proof}
Let $\del: \mathrm{Hom}_S(-,X) \sra \calS_X$ denote the map of functors induced from the diagonal $\Delta \in \calS_X(X)$, and denote by $p_1$ the first projection $X \times_S X \sra X$.

Then for any $S$-morphism $g: T \sra X$ we have $\del(g) = (p_1 \circ \Delta \circ g, 1_T) = (g, 1_T)$, as seen by the diagram below: 
\[
\xymatrix{
X \times_S T \ar[r] \ar[d] & X \times_S X \ar[r]^-{p_1} \ar[d] & X \ar[d]\\
T \ar[r]^g \ar@/^1pc/[u]^{(g,1_T)} & X \ar[r] \ar@/_1pc/[u]_{\Delta} & S. 
}
\]
Also, for any section $\xi = (g,1_T) : T \sra X \times_S T$ we have by definition $n(\xi) = g$. It is thus clear that $n$ and $\del$ are inverse to each other.
\end{proof}

\begin{rem}
It follows from Lemma (\ref{lem:sec_to_hom}) that any section of a separated morphism $f : X\sra S$ is closed, being the pull-back of the diagonal $\Delta: X\sra X\times_SX$.
\end{rem}

\begin{lemma}
\label{lem:infinitesimal_sections}
Let $f : X \sra \Spec (A)$ be a morphism of schemes, where  $(A,\m)$ is a complete local ring. For each $n\geq 1$ we let $A_n:= A/\m^{n+1}$ and $X_n:=X\times_A \Spec(A_n)$. Furthermore we let $f_n : X_n  \sra \Spec (A_n)$ denote the induced morphism.
\begin{itemize}
\item[(a)] Sections $\xi_n$ of $f_n : X_n \sra \Spec (A_n)$ are closed.
\item[(b)] Let $\xi '$ and $\xi$ be sections of $f : X \sra \Spec (A)$,  such that when restricted to $\Spec (A_n)$ we have $\xi_n'=\xi_n$ for all $n$. Then $\xi=\xi'$.
\end{itemize}
\end{lemma}

\begin{proof}
As the underlying topological space of  $\Spec (A_n) $ is one point it follows that $\xi_n$ factors through any open affine $U\subseteq X_n$ containing the image of $\xi_n$. Then $\Spec (A_n) \sra U$ is a section of the separated morphism $f_n | U$ and hence a closed section. It follows that that the image of $\xi_n$ is closed in any open affine $U\subseteq X_n$, and hence closed in $X_n$. This proves the first statement.

To prove the second assertion we let $\Spec (B)\subseteq X$ be an open affine subscheme containing the image of the closed point $\Spec (A/\m)$ under $\xi $. The sections $\xi_n : \Spec (A_n) \sra X_n$, composed with the closed immersions $X_n \sra X$ factor through $\Spec (B)$. The corresponding ring homomorphism $B\sra A_n$ determines a unique morphism to the inverse limit $\varprojlim A_n=A$. As $\xi '$ coincides with $\xi$ when restricted to $\Spec (A_n)$ it follows that $\xi'=\xi$.
\end{proof}

\subsection{The Hilbert functor of points}

For a fixed scheme $X$ over some base $S$ we let $\calH^m_X$ be the Hilbert functor of $m$-points on $X$, as defined by Grothendieck \cite[p. 221-26]{fga}. Thus, for any $S$-scheme $T$ we have that the $T$-valued points of $\calH^m_X$ is the set of closed subschemes $Z\subseteq X \times_S T$ such that the induced projection $p:Z\sra T$ is flat and finite of rank $m$. In other words, $p$ is flat and finite, and $p_* \calO_Z$ is locally free of rank $m$ as an $\calO_T$-module.

\begin{lemma}
Let $U\subseteq X$ be an open  subscheme, where $X$ is separated over the base scheme $S$. Then there is a natural transformation $\calH_U^m \sra \calH_X^m$.
\end{lemma}

\begin{proof}
As a $T$-valued point $Z$ of $\calH_U^m$ is finite over $T$ it is in particular proper over $T$ by \cite[Cor. 6.1.11]{egaII}. Thus the composition 
\[
Z \ra U \times_S T \ra X \times_S T \ra T
\]
is proper and the projection map $X \times_S T \sra T$ is separated. It follows from \cite[Cor. 5.4.3]{egaII} that the immersion $Z \sra X \times_S T$ is proper and hence closed, and so we get an element of $\calH_X^m(T)$.
\end{proof}

\begin{rem}
\label{rem:no_map}
If $U\subset X$ is an open immersion then it is \emph{not} true in general that we have a map of functors $\calH^1_{U}\sra \calH^1_{X}$, as the following example shows. 
\end{rem}

\begin{exam}
Let $X$ denote the line with a double point. We obtain $X$ by glueing two copies of the line  along the open complement of a closed point. In particular we can let $U\subset X$ be the open subscheme given by one of the lines. Finally we let the base $S$ be the line, with the natural projection morphism $X \sra S$. Now it is clear that the whole line $U$ itself is flat and finite of rank one over the base. However, the line $U$ is not closed in $X$, but clearly closed in $U$. Thus there is not a natural map from $\calH^1_U $ to $\calH^1_X$. In fact $\calH_U^1(U)$ is a singleton set, whereas $\calH_X^1(U)$ is the empty set.
\end{exam}

\begin{rem}
Let $\calF$ be a coherent sheaf on $X$, and let $\mathrm{Quot}(\calF/X/S)$ denote the Quot-functor (see e.g. \cite{fga} or \cite{artin_alg_form_moduli_I}). When $f : X \sra S$ is locally of finite presentation, Artin applies the algebraization theorem to show that the Quot-functor is representable by an algebraic space in \cite[\S 6]{artin_alg_form_moduli_I} - but the additional hypothesis that $f : X \sra S$ is separated is required. 

One instance where separatedness is needed in the proof of \cite[Thm 6.1, p.61]{artin_alg_form_moduli_I} is a reduction to the Quot-functor $\mathrm{Quot}(\calF'/X'/S)$, where $X'$ is open in $X$. Because, as pointed out in Remark (\ref{rem:no_map}), there does not always exist a map $\mathrm{Quot}(\calF'/X'/S) \sra \mathrm{Quot}(\calF/X/S)$.

Also, the reference to Grothendieck's existence theorem for formal sheaves \cite[Thm. 5.1.4]{egaIII} in \cite[Thm 6.1, p.64]{artin_alg_form_moduli_I} is problematic, since the existence theorem requires separatedness conditions.

The problems of separatedness are adressed by Artin himself in \cite[Appendix, p. 186]{artin_versal_def}.
\end{rem}

\subsection{The Hilbert functor of one point}
\label{sec:hilb_one}

We will now focus on a particular case of the Hilbert functor of points, namely when $m = 1$. In that case we have that the projection $p: Z \sra T$ is finite and flat of rank 1, and then $p$ must  be an isomorphism. The inverse to $p$ gives a closed section $T \sra X \times_S T$ and thus we may identify the set $\calH^1_X(T)$ with the set of closed sections $T \sra X \times_S T$. In particular we see that $\calH^1_X$ is a subfunctor of the section functor  $\calS_X$ (\ref{sec:sections}).

\begin{prop}
\label{prop:rep_hilb1}
The map of functors (\ref{eq:sec_to_hom}) induces a natural map
\eqbeg
\label{eq:norm_map}
n_X : \calH_X^1 \ra \mathrm{Hom}_{S}(-,X),
\eqend
which is an isomorphism if and only if $f : X\sra S$ is separated.
\end{prop}

\begin{proof}
When $f : X\sra S$ is separated we have that any section is closed. Consequently we have that  $\calH_X^1=\calS_X$, and the proposition is then a special case of Lemma (\ref{lem:sec_to_hom}). When $f : X\sra S$ is not separated there exists non-closed sections, e.g. the diagonal map $X \sra X\times_S X$. Therefore $\calH_X^1 (X) \subset \calS_X(X)$ is a proper subset and the map $n_X$ is not an isomorphism.
\end{proof}

\begin{rem}
In \cite[p. 221-26]{fga} Grothendieck introduced a norm map from the Hilbert scheme $\mathrm{Hilb}^m_{X}$ to the $m$-fold symmetric product $\mathrm{Sym}^m_S(X)$. The map (\ref{eq:norm_map}) is this norm map for $m=1$.
\end{rem}

\begin{rem}
In the definition of the Hilbert functor $\calH^m_X$, one could replace closed subschemes with locally closed subschemes. In that case we would have equality $\calH_X^1=\calS_X$, and in particular the norm map (\ref{eq:norm_map}) would be an isomorphism. See \cite[Prop 2.2, Cor. 2.3]{kleiman_mult_point_II}. However, it is not clear that the refined definition of $\calH_X^m$ would prove to be representable (for $m>1$). 

For other discussions of Hilbert- and Quot functors related to stacks see \cite[p. 186]{artin_versal_def}, \cite{olsson_starr} and  \cite{vistoli_hilbert_stack}.
\end{rem}

\section{Formal deformations of the Hilbert functor}

Before we prove our main result, we will, for the sake of completeness, prove that deformations of algebraic spaces are effective. A proof of this can also be found in \cite{artin_alg_form_moduli_I}.

\subsection{Effective deformations}

Let $S$ be a scheme, and let $F$ be a contravariant functor from the category of $S$-schemes to sets. If $X = \Spec (A)$ is an affine scheme, we write $F(A)$ instead of $F(X)$.

Given a field $k$ and an element $\xi_0 \in F(k)$. A \emph{formal deformation} of $\xi_0$ is a pair $(A, \{\xi_n\}_{n \geq 0})$ where $A$ is a complete local ring with residue field $k$ and $\{\xi_n\}_{n \geq 0}$ is a collection of elements with $\xi_n \in F(A / \m^{n+1})$ such that $\xi_{n-1}$ is induced from $\xi_{n}$ and $\xi_0$ is the original element. The deformation is called \emph{effective} if there is an element $\xi \in F(A)$ inducing the elements $\{ \xi_n \}$.

\begin{rem}
\label{rem:formal_defs_of_schemes}
Any formal deformation of a scheme $F=\mathrm{Hom}_S(-,Y)$ is effective. Indeed, if $\xi_n : \Spec (A/\m^{n+1})\sra Y$ is a compatible collection of morphisms, then all $\xi_n$ factor through any open affine $U\subseteq Y$ containing the image of the point $\Spec (A/\m)$. Consequently the collection of maps $\xi_n$ can be reduced to the affine case where the result follows from the universal property of the inverse limit $\varprojlim A/\m^{n+1}$.
\end{rem}

\begin{prop}
Any formal deformation of an algebraic space $X$ is effective.
\end{prop}

\begin{proof}
We can carry over the arguments of Remark \ref{rem:formal_defs_of_schemes} to the setting of algebraic spaces. Thus suppose that we are given a collection $\mas{\xi_n}$ of morphisms $\xi_n: \Spec(A / \m^{n+1}) \sra X$. By \cite[Thm. II.6.4]{knutson_alg_spaces} we have that the map $\xi_0: \Spec(k) \sra X$ factors as $\Spec(k) \sra U \sra X$ where $U$ is affine and $U \sra X$ is étale.

By the lifting property of étale maps we have that also the maps $\xi_n: \Spec(A / \m^{n+1}) \sra X$ factor as $\Spec(A / \m^{n+1}) \sra U \sra X$. Since $U$ is affine it is clear that we have a map $\Spec(A) \sra U$ that restricts to the given maps $\Spec(A / \m^{n+1}) \sra U$ for each $n$, and so we have found the required map $\xi: \Spec(A) \sra X$. Thus the deformation is effective.
\end{proof}

\subsection{Formal deformations of the Hilbert functor}

For any complete local ring $(A,\m)$ we consider the natural map
\eqbeg
\label{eq:invlim_map}
\calH^1_X(A) \ra \varprojlim \calH_X^1(A/\m^{n+1}), 
\eqend
where $\calH_X^1$ is the Hilbert functor defined in (\ref{sec:hilb_one}). A consequence of Lemma (\ref{lem:infinitesimal_sections}) is that the map (\ref{eq:invlim_map}) is injective. Our main result is that that the map (\ref{eq:invlim_map}) is not always surjective.

\subsection{Valuative criterion of separatedness} Recall that a morphism of schemes $f : X \sra S$ satisfies the valuative criterion for separatedness if, for any valuation ring $A$ with fraction field $K$, and for any commutative diagram of schemes

\eqbeg
\label{eq:val_criterion}
\xymatrix{
\Spec (K) \ar[r] \ar[d]^{\can} & X \ar[d]^{f} \\
\Spec (A) \ar[r]   & S,
}
\eqend
there exists at most one morphism $ \xi : \Spec (A) \sra X$ extending the morphism $\Spec (K) \to X$ of the diagram.  
\begin{thm}
\label{thm:hilbone_not_rep} Let $f : X\sra S$ be a  morphism locally of finite type, and where the base scheme $S$ is locally noetherian. Assume furthermore that $f : X\sra S$  does not satisfy the valuative criterion of separatedness. Then the Hilbert functor $\calH_X^1$ has non-effective formal deformations. In particular, for such $X$, the functor $\calH_X^1$ is not representable by a scheme or an algebraic space.
\end{thm}

\begin{proof} As $f : X \sra S$ does {\em not} satisfy the valuative criterion of separatedness there exists a diagram as (\ref{eq:val_criterion}) with at least two extensions $\Spec(A) \to X$ of the morphism from the generic point $\Spec (K) \to X$. Furthemore, as $f : X \sra S$ is locally of finite type and $S$ is locally noetherian, we can assume that $A$ is a discrete valuation ring \cite[Prop. 7.2.3]{egaII}. We have that the completion $\hat{A}$ of $A$ is also a discrete valuation ring with fraction field $\hat{K}$ and the map $A \sra \hat{A}$ is injective. Through the canonical map $\Spec(\hat{A}) \sra \Spec(A)$ we then obtain a diagram of the form (\ref{eq:val_criterion}) with $\hat{A}$ and $\hat{K}$ instead of $A$ and $K$. Thus we may assume that the ring $A$ is complete.
 
Consider now one of the sections $\xi : \Spec (A) \sra X_A:=X\times_S \Spec (A)$ that we obtain from the diagram (\ref{eq:val_criterion}). As we have at least one other section extending the map from the generic point $\Spec (K)\sra X_A$ it follows that $\xi $ is not closed. Hence $\xi$ is not a $\Spec (A)$-valued point of $\calH_X^1$.

By Lemma (\ref{lem:infinitesimal_sections}), the induced sections $\xi_n : \Spec (A/\m^{n+1}) \sra X\times _S \Spec (A/\m^{n+1})$ are closed. Therefore the collection of sections $\{\xi_n\}$ form a formal deformation of $\xi_0 : \Spec (A/\m) \sra X\times_S \Spec (A/\m)$. If there was an element $\xi ' \in \calH^1_X(A)$ whose truncations would form the constructed formal deformation, it would have to be $\xi$ by Lemma (\ref{lem:infinitesimal_sections}). However, as $\xi \notin \calH^1_X(A)$ we have that the constructed deformation is not effective.
\end{proof}

\begin{cor} Let $f : X \sra S$ be locally of finite type, with $S$ locally noetherian. Then $\calH^1_X$ is representable if and only if $f : X \to S$ is separated. If it exists, the representing scheme is $X$.
\end{cor}
\begin{proof} By \cite[Prop. 7.2.3, Rem. 7.2.4(i)]{egaII} we have that $f : X \sra S$ is separated if and only if the valuative criterion holds. Consequently the corollary follows from Proposition (\ref{prop:rep_hilb1}) and Theorem (\ref{thm:hilbone_not_rep}).
\end{proof}

\begin{rem} When $f : X\sra S$ is locally of finite type of locally noetherian schemes, the valuative criterion (\ref{eq:val_criterion}) can be checked using complete discrete valuation rings $A$ \cite[Prop.7.2.3]{egaII}. For general $ f : X\sra S$ the valuative criterion can be checked using ``complete'' valuation rings $A$ (see \cite[Remark. 7.2.4 (ii)]{egaII}). However, the terminology ``complete'' in \cite{egaII} does not mean that the ring $A$ is Hausdorff and complete in the $\m$-adic topology, unless $A$ is a DVR. In particular, for these rings the canonical map $A \sra \varprojlim A/ \m^{n+1}$ need not be an isomorphism which is the defining property of completeness used in this article. We thank David Rydh who made us aware of this.
\end{rem}


\bibliographystyle{dary}

\bibliography{papers}

\end{document}